\documentclass[11pt]{amsart}

\usepackage[T1]{fontenc}
\usepackage{lmodern}
\usepackage{amsmath,amssymb,amsthm,mathtools}
\usepackage{microtype}
\usepackage[colorlinks=true,linkcolor=blue,citecolor=blue,urlcolor=blue]{hyperref}
\hypersetup{
  pdftitle={Non-Unique Factorization in Polynomial Rings over Division Rings},
  pdfauthor={Elad Paran},
  pdfsubject={Unique factorization, half-factoriality, and polynomial identities over division rings}
}

\newtheorem{theorem}{Theorem}[section]
\newtheorem{proposition}[theorem]{Proposition}
\newtheorem{lemma}[theorem]{Lemma}

\theoremstyle{definition}

\theoremstyle{remark}

\newcommand{\Cent}{\operatorname{C}}

\title{Non-Unique Factorization in Polynomial Rings over Division Rings}
\author{Elad Paran}
\address{Department of Mathematics and Computer Science, The Open University of Israel, Ra\'anana 4353701, Israel}
\email{paran@openu.ac.il}
\date{}

\subjclass[2020]{Primary 16U30; Secondary 16R20}
\keywords{division ring, polynomial ring, similarity, unique factorization, half-factoriality, polynomial identity}

\begin{document}

\begin{abstract}
Let $D$ be a division ring and let $R=D[x_1,\ldots,x_r]$ be the ring of polynomials in $r$ central variables over $D$.  We prove that, for $r\geq 2$, $R$ is a unique factorization domain in the sense of Cohn if and only if $D$ is commutative.  Moreover, if $D$ is noncommutative, then $R$ is not half-factorial.  This generalizes the rational-quaternion case proved by Beauregard and the real-quaternion case proved by Skopenkov and Krasauskas.
\end{abstract}

\maketitle

\section{Introduction}

Let $D$ be a division ring and let $R = D[x_1,\ldots,x_r]$ be the ring of polynomials in $r$ central variables over $D$. If $D$ is a field, then $R$ is of course a unique factorization domain (UFD). If $D$ is not a field and $r = 1$, then $D[x_1]$ is still a UFD in Cohn's sense: factorizations into atoms (irreducible nonunits) are unique up to order and similarity (the precise definition of similarity is recalled below). This is a classical result of Ore  \cite{Ore}. All polynomial rings considered below are atomic, by induction on total degree. For $r > 1$ the situation turns out to be different. A fundamental counterexample is due to Beauregard, who proved that $D[x,y]$ is not a UFD when $D$ is the rational Hamilton quaternion algebra \cite{Beauregard}.  His proof is based on the polynomial
\[
 f=(x^2y^2-1)+(x^2-y^2)i+2xyj,
\]
which is irreducible over the rational quaternions and for which $f\overline f$ admits atomic factorizations of different lengths. This means that in this case, the ring $D[x,y]$ is not even {\it half-factorial} (an atomic ring is called half-factorial if every factorization of a given element into atoms has the same length). The irreducibility of the same polynomial over the real Hamilton quaternions was proved later by Skopenkov and Krasauskas \cite{SkopenkovKrasauskas}, which showed the same failure of half-factoriality in the real case. Several other works have studied aspects of factorization for specific cases of polynomial rings over division rings; see, for example, \cite{FrischaufLubbesSchroecker,GoutorTikhonov,Grzegorzewski,Koprowski,LercherEtAl,LercherSchroecker,LiScharlerSchroecker}.

These examples raise the following question: When is $D[x_1,\ldots,x_r]$ a UFD? To the author's knowledge, this question has not previously been answered in this generality. Here we give a complete answer:

\begin{theorem}\label{thm:main}
Let $D$ be a division ring and let $r\geq 2$. Then $D[x_1,\ldots,x_r]$ is a UFD if and only if $D$ is commutative.

\end{theorem}

In fact, we prove a sharper statement:

\begin{theorem}\label{thm:half-factorial}
Let $D$ be a noncommutative division ring, and let $r\geq 2$.  Then $D[x_1,\ldots,x_r]$ is not half-factorial.
\end{theorem}

This paper is organized as follows.  In Section~2 we prove Theorem~\ref{thm:main} by constructing, from two noncommuting elements of $D$, two inequivalent factorizations into atoms of the form $x-H(y)$; the key observation is that similarity of such atoms is equivalent to conjugacy of the corresponding polynomials $H(y)$.  In Section~3 we begin the proof of Theorem~\ref{thm:half-factorial} with a commutator construction producing factorizations of lengths two and three, and reduce the remaining case, by a polynomial identity and Kaplansky's theorem, to quaternion division algebras.  In Section~4 we complete the proof by explicit norm constructions, treating separately characteristic two and characteristic different from two.

\section{Proof of Theorem \ref{thm:main}}

Recall that two nonzero elements $a,b$ of a domain $R$ are \emph{similar} if
\[
 R/aR\cong R/bR
\]
as right $R$-modules.  We say that $R$ is a \emph{unique factorization domain}, or a \emph{UFD}, if it is atomic and every equality
\[
 p_1\cdots p_m=q_1\cdots q_n
\]
between products of atoms has $m=n$, and, after a permutation, $p_i$ is similar to $q_i$ for every $i$.  We refer to \cite{CohnUFD,CohnBook,Smertnig} for background.

We use two observations, stated with several coefficient variables for later application to $D[x_1,\ldots,x_r]$.

\begin{lemma}\label{lem:atom}
Let
\[
 S=D[y_1,\ldots,y_s],\qquad R=S[x],
\]
where all variables are central, and let $H\in S$.  Then $x-H$ is an atom of $R$.
\end{lemma}

\begin{proof}
Write $x-H=AB$.  By additivity of $x$-degree, at least one of $A,B$ belongs to $S$.  If $A\in S$, comparison of leading coefficients in $x$ makes $A$ right-invertible in $S$, which implies that $A\in D^\times$. The case $B\in S$ is symmetric.
\end{proof}

\begin{lemma}\label{lem:similarity}
With the notation of Lemma~\ref{lem:atom}, let $H,K\in S$.  Then the following conditions are equivalent:
\begin{enumerate}
\item $x-H$ is similar to $x-K$;
\item $K=hHh^{-1}$ for some $h\in D^\times$.
\end{enumerate}
\end{lemma}

\begin{proof}
Right division by $x-H$ gives
\[
 R/(x-H)R\cong S,
\]
as right $S$-modules, with $x$ acting by $s\mapsto Hs$.  An isomorphism between the two quotient modules is therefore left multiplication by some $h\in S^\times=D^\times$, and it is compatible with the action of $x$ precisely when
\[
 hH=Kh,
\]
as required.
\end{proof}

\begin{theorem}\label{thm:universal-construction}
If $D$ is noncommutative, then $D[x,y]$ is not a UFD.
\end{theorem}

\begin{proof}
Choose $c,v\in D$ such that $cv\neq vc$.  In $D[y]$, put
\[
 U=y+v,\qquad
 P=Uc+vy^2,\qquad
 Q=cU+vy^2,\qquad
 T=y^3U.
\]
Since $U$ commutes with $vy^2$, we have $UQ=PU$, and hence
\[
 TQ=y^3UQ=y^3PU=PT.
\]
Expanding now gives
\begin{equation}\label{eq:universal}
 (x-P-T)(x-Q)=(x-P)(x-Q-T).
\end{equation}
All four factors are atoms by Lemma~\ref{lem:atom}.  Moreover,
\[
 \deg_y P=\deg_y Q=2,
 \qquad
 \deg_y(P+T)=\deg_y(Q+T)=4.
\]
Thus Lemma~\ref{lem:similarity} forces the pairings $x-P\sim x-Q$ and $x-P-T\sim x-Q-T$.  But if $Q=hPh^{-1}$ for some $h\in D^\times$, then comparison of the coefficients of $y^2$, $y$, and $1$ gives 
\[
 hv=vh,\qquad hc=ch,\qquad vc=cv,
\]
a contradiction.  Hence the factorizations in \eqref{eq:universal} are inequivalent.
\end{proof}

\begin{proof}[Proof of Theorem~\ref{thm:main}]
For commutative $D$ this is the usual polynomial-ring theorem.  If $D$ is noncommutative, apply Theorem~\ref{thm:universal-construction} to $D[x_1,x_2]$.  Lemmas~\ref{lem:atom} and \ref{lem:similarity}, with $S=D[x_2,\ldots,x_r]$ and $x=x_1$, preserve both atomicity and the similarity obstruction in $D[x_1,\ldots,x_r]$.
\end{proof}

\section{A commutator construction and a PI reduction}

The next construction is valid over an arbitrary division ring.

For $a\in D$, we denote by $\Cent_D(a)$ the centralizer of $a$ in $D$.

\begin{lemma}\label{lem:quadratic-atom}
Let $a,b\in D$ and suppose that $c=[a,b]\neq 0$.  Put
\[
 F=(y+b)c^{-1}(x+a)-1.
\]
Then $F$ is an atom of $D[x,y]$ if and only if
\[
 cbc^{-1}+b\notin \Cent_D(a).
\]
\end{lemma}

\begin{proof}
We have
\begin{equation}\label{eq:F-expand}
 F=c^{-1}xy+bc^{-1}x+c^{-1}ay+bc^{-1}a-1.
\end{equation}
In a proper factorization both factors have degree one.  Their homogeneous parts, in one order or the other, are an $x$-term and a $y$-term.  For
\[
 F=(uy+r)(px+s),
 \qquad u,p\neq 0,
\]
coefficient comparison gives
\[
 up=c^{-1},\qquad us=c^{-1}a,
 \qquad rp=bc^{-1}.
\]
Hence $rs=bc^{-1}a$, a contradiction.  For
\[
 F=(ux+r)(qy+s),
 \qquad u,q\neq 0,
\]
we obtain
\[
 uq=c^{-1},\qquad us=bc^{-1},
 \qquad rq=c^{-1}a,\qquad
 rs=c^{-1}acbc^{-1}.
\]
Conversely, these equations are realized by
\[
 u=1,\qquad q=c^{-1},\qquad r=c^{-1}ac,
 \qquad s=bc^{-1}.
\]
Thus $F$ is reducible precisely when
\[
 c^{-1}acbc^{-1}=bc^{-1}a-1.
\]
Multiplication by $c$ on the left, together with $c=[a,b]$, makes this equivalent to
\[
 [a,cbc^{-1}+b]=0.
\]
\end{proof}

\begin{proposition}\label{prop:commutator-obstruction}
Suppose that there are $a,b\in D$ such that, for $c=[a,b]$, one has
\[
 c\neq 0
 \qquad\text{and}\qquad
 cbc^{-1}+b\notin \Cent_D(a).
\]
Then $D[x,y]$ is not half-factorial.
\end{proposition}

\begin{proof}
Let $F$ be the atom in Lemma~\ref{lem:quadratic-atom}, and put
\[
 A=x+a,\qquad B=y+b.
\]
Since $AB=BA+c$, we have
\[
 (x+a)F=(y+b)\bigl((x+a)c^{-1}\bigr)(x+a).
\]
The factor $F$ is an atom by Lemma~\ref{lem:quadratic-atom}; the other factors have degree one.  These are atomic factorizations of lengths two and three.
\end{proof}

\begin{lemma}\label{lem:pi-reduction}
Suppose that
\begin{equation}\label{eq:exceptional-case}
 cbc^{-1}+b\in \Cent_D(a)
\end{equation}
for all $a,b\in D$ for which $c=[a,b]\neq 0$.  Then $D$ is finite-dimensional over its center.  If $D$ is noncommutative, its central degree is two; in particular, $D$ is a quaternion division algebra over its center.
\end{lemma}

\begin{proof}
For $c=[a,b]\neq 0$, multiplying $[a,cbc^{-1}+b]=0$ on the left by $c^{-1}$ and on the right by $c$, and using $c=[a,b]$, gives
\begin{equation}\label{eq:exceptional-symmetric}
 c^{-1}ac+a\in \Cent_D(b).
\end{equation}
Put $q=c^{-1}ac$ and $s=a+q$.  Equation \eqref{eq:exceptional-symmetric} says that $s$ commutes with $b$.  Since
\[
 [a,ab]=ac\neq 0,
\]
applying it to $(a,ab)$ gives $[s,ab]=0$, because
\[
 (ac)^{-1}a(ac)=c^{-1}ac=q;
\]
therefore
\[
 0=[s,ab]=[s,a]b,
\]
so $[s,a]=0$.  Hence $[s,c]=0$, and $sc=cs$ gives
\[
 ac^2=c^2a.
\]
Thus, also when $[a,b]=0$,
\[
 [[a,b]^2,a]=0
 \qquad(a,b\in D).
\]
The polynomial $[[X,Y]^2,X]$ is nonzero over every prime field: its reduced expansion has six distinct monomials.  Since $D$ is primitive, Kaplansky's theorem \cite[Theorem~1]{KaplanskyPI} implies that $D$ is finite-dimensional over its center.

Assume now that $D$ is noncommutative.  Let $K$ be its center, and write $[D:K]=n^2$.  The field $K$ is infinite, for otherwise $D$ would be finite and hence commutative.  Consequently every polynomial identity of $D$ remains an identity after scalar extension: this follows by expanding the variables in a $K$-basis of $D$ and comparing the resulting polynomial coefficients over the infinite field $K$.  Over an algebraic closure $\overline K$ we have
\[
 D\otimes_K\overline K\cong M_n(\overline K).
\]
If $n\geq 3$, take, in the upper-left $3\times3$ block,
\[
 X=E_{23}+E_{32},\qquad Y=E_{31}+E_{33}.
\]
Then
\[
 [X,Y]=E_{21}+E_{23}-E_{32},\qquad
 [X,Y]^2=-E_{22}-E_{31}-E_{33},
\]
and hence
\[
[[X,Y]^2,X]=E_{21}\neq 0.
\]
Thus $n\leq 2$.  Since $D$ is noncommutative, $n=2$, and $D$ is a quaternion division algebra over $K$.
\end{proof}

\section{Quaternion division algebras}

Central degree two requires separate constructions in the two possible characteristics.

\begin{proposition}\label{prop:quaternion-odd}
Let $K$ have characteristic different from two, and let
\[
 D=(a,b)_K=K+Ki+Kj+Kij
\]
be a quaternion division algebra, where
\[
 i^2=a,\qquad j^2=b,\qquad ji=-ij.
\]
Then $D[x,y]$ is not half-factorial.
\end{proposition}

\begin{proof}
Put
\begin{equation}\label{eq:quat-odd-F}
 F=x^2y^2+a^2b+(bx^2-ay^2)i+2axyj.
\end{equation}
Writing
\[
 A=x^2y^2+a^2b,\qquad B=bx^2-ay^2,
 \qquad C=2axy,
\]
we obtain
\begin{align}
 F\overline F
 &=A^2-aB^2-bC^2 \notag\\
 &=(x^4-a^3)(y^4-ab^2) \notag\\
 &=(x^2-ai)(x^2+ai)(y^2-bi)(y^2+bi).
 \label{eq:quat-odd-norm}
\end{align}
If a factorization had a factor independent of $y$, that factor would be a nonconstant left or right divisor of both
\[
 2axj\qquad\text{and}\qquad x^2-ai.
\]
This is impossible, since a nonconstant divisor of $2axj$ has zero constant coefficient.  The same argument, with $x$ and $y$ interchanged, excludes every one-variable factor.  Any proper factorization therefore has, after transferring a constant unit from one factor to the other, the form
\[
 F=(xy+px+qy+r)(xy-px-qy+u).
\]
Its coefficients give
\begin{gather}
 p^2=-bi,\qquad q^2=ai, \label{eq:odd-pq}\\
 pu=rp,\qquad qu=rq, \label{eq:odd-intertwine}\\
 u+r-pq-qp=2aj. \label{eq:odd-xy}
\end{gather}
Now $p,q\in\Cent_D(i)=K(i)$.  Put $h=q^{-1}p$.  By \eqref{eq:odd-intertwine}, $hu=uh$.  If $h\notin K$, then $\Cent_D(h)=K(i)$, so $u,r\in K(i)$, contrary to \eqref{eq:odd-xy}.  If $h=d\in K^\times$, then $p=qd$ and
\[
 d^2=-\frac ba.
\]
Hence
\[
 (i+d^{-1}j)^2=a+d^{-2}b=0,
\]
a contradiction.  Thus $F$ is an atom.  Since the canonical involution reverses products, $\overline F$ is an atom as well.  Equation \eqref{eq:quat-odd-norm} gives an atomic factorization of length two and, after factoring its four nonunit factors, one of length at least four.
\end{proof}

\begin{proposition}\label{prop:quaternion-two}
Let $K$ have characteristic two, and let
\[
 D=[a,b)_K
\]
be a quaternion division algebra.  Here $[a,b)_K$ denotes the
characteristic-two quaternion symbol algebra
\[
 D=K+Ki+Kj+Kij,
\]
where
\[
 i^2+i=a,\qquad j^2=b,\qquad ji=(i+1)j.
\]
Then $D[x,y]$ is not half-factorial.
\end{proposition}

\begin{proof}
Set
\begin{equation}\label{eq:quat-two-F}
 F=x^2y^2+ab+(bx^2+y^2)i+xyj.
\end{equation}
The canonical involution satisfies $\overline i=i+1$ and $\overline j=j$.  Writing
\[
 A=x^2y^2+ab,\qquad B=bx^2+y^2,
 \qquad C=xy,
\]
gives
\begin{align}
 F\overline F
 &=A^2+AB+aB^2+bC^2 \notag\\
 &=(x^4+x^2+a)(y^4+by^2+ab^2) \notag\\
 &=(x^2+i)(x^2+i+1)(y^2+bi)(y^2+b(i+1)).
 \label{eq:quat-two-norm}
\end{align}
If a factorization had a factor independent of $y$, that factor would be a nonconstant left or right divisor of both $xj$ and $x^2+i$.  A nonconstant divisor of $xj$ has zero constant coefficient, so it cannot divide $x^2+i$.  Interchanging $x$ and $y$ and using the pair $yj,y^2+bi$ excludes a factor independent of $x$.  Thus any proper factorization has, after transferring a constant unit from one factor to the other, the form
\[
 F=(xy+px+qy+r)(xy+px+qy+u).
\]
Its coefficients give
\begin{gather}
 p^2=bi,\qquad q^2=i, \label{eq:two-pq}\\
 pu=rp,\qquad qu=rq, \label{eq:two-intertwine}\\
 u+r+pq+qp=j. \label{eq:two-xy}
\end{gather}
Now $p,q\in K(i)$.  Put $h=q^{-1}p$.  If $h\notin K$, then $\Cent_D(h)=K(i)$, so \eqref{eq:two-intertwine} gives $u,r\in K(i)$, contrary to \eqref{eq:two-xy}.  If $h=d\in K^\times$, then $d^2=b$ and
\[
 (d^{-1}j+1)^2=0,
\]
a contradiction.  Hence $F$ is an atom.  Since the canonical involution reverses products, $\overline F$ is an atom as well, and \eqref{eq:quat-two-norm} gives atomic factorizations of lengths two and at least four.
\end{proof}

\begin{proof}[Proof of Theorem~\ref{thm:half-factorial}]
If the hypothesis of Proposition~\ref{prop:commutator-obstruction} holds, we are done.  Otherwise the condition of Lemma~\ref{lem:pi-reduction} holds for every noncommuting pair.  Hence $D$ is a quaternion algebra, and Propositions~\ref{prop:quaternion-odd} and \ref{prop:quaternion-two} apply.

For $r>2$, degree in each added variable shows that the atoms used in $D[x_1,x_2]$ remain atoms in $D[x_1,\ldots,x_r]$.
\end{proof}

\section*{Declarations}

\subsection*{Declaration of AI use}
This paper was written with AI assistance, model ChatGPT 5.6-Sol. The central result, Theorem 1.1 was first obtained by the author, with the same proof given here. This proof in fact precedes the current AI era. This led to a short note written by the author, that has not been submitted for publication. The note was uploaded for review to an AI model in September 2026. The model was able to sharpen the author's result and prove Theorem 1.2 concerning half-factoriality. The model's proof was then verified by the author, who takes full responsibility for the entire contents of the paper.

\subsection*{Data availability}

Data sharing is not applicable to this article as no datasets were generated or analyzed during the current study.

\subsection*{Competing interests}

The author has no relevant financial or non-financial interests to disclose.

\end{document}